\documentclass[11pt, a4paper]{article}
\usepackage{amsmath, amsthm, amssymb, url, longtable, changes, enumitem}

\theoremstyle{plain}
\newtheorem{theorem}{Theorem}[section]

\newtheorem{corollary}[theorem]{Corollary}
\newtheorem{lemma}[theorem]{Lemma}

\newtheorem*{claim*}{Claim}
\newtheorem*{problem*}{Problem}
\newtheorem*{conjecture*}{Conjecture}

\theoremstyle{definition}
\newtheorem{definition}[theorem]{Definition}

\newcommand{\A}{\mathrm{A}}
\newcommand{\B}{\mathrm{B}}
\newcommand{\C}{\mathrm{C}}
\newcommand{\la}{\langle}
\newcommand{\ra}{\rangle}

\DeclareMathOperator{\Ima}{Im}
\setlist[enumerate,1]{label={\upshape\arabic*.}}

\usepackage{ltablex} % So tabularx tables can break over pages
\usepackage{tabularx}
\newcolumntype{C}[1]{>{\centering\arraybackslash}m{#1}}
\title{$3$-generated axial algebras with a minimal Miyamoto group}
\author{J.~M\textsuperscript{c}Inroy\footnote{School of Mathematics, University of Bristol, Fry Building, Woodland Road, Bristol, BS8 1UG, UK; and the Heilbronn Institute 
for Mathematical Research, Bristol, UK; email: 
justin.mcinroy@bristol.ac.uk}}

\date{\today}
\begin{document}
\maketitle

\begin{abstract}
Axial algebras are a recently introduced class of non-associative algebra, with a naturally associated group, which generalise the Griess algebra and some key features of the moonshine VOA.  Sakuma's Theorem classifies the eight $2$-generated axial algebras of Monster type.  In this paper, we compute almost all the $3$-generated axial algebras whose associated Miyamoto group is minimal $3$-generated (this includes the minimal $3$-generated algebras).  We note that this work was carried out independently to that of Mamontov, Staroletov and Whybrow and extends their result by computing more algebras and not assuming primitivity, or an associating bilinear form.
\end{abstract}

\section{Introduction}

There has been much interest recently in axial algebras and in particular their construction.  There are eight $2$-generated axial algebras of Monster type and these are classified in Sakuma's Theorem \cite{Axial1}.  In this note, we look at $3$-generated axial algebras.

Axial algebras are non-associative algebras which axiomatise some key features of VOAs and the Griess algebra.  For full details, we refer readers to \cite{Axial1}, or for a more general treatment \cite{axialstructure}.  Roughly speaking, they are generated by a set $X$ of distinguished elements, called \emph{axes}.  The multiplicative action of an axis on the algebra decomposes it into a direct sum of eigenspaces and the multiplication of these eigenspaces is governed by a so-called fusion law.  This gives us partial control over multiplication in the algebra.  Importantly, when this fusion law is graded, we can associate automorphisms, called \emph{Miyamoto automorphisms}, to each axis.  The group generated by all of these is called the \emph{Miyamoto group}.

The Griess algebra is the prototypical example of an axial algebra and its fusion law is given in Table \ref{tab:monsterfusion}.  However, there are many other axial algebras with this fusion law, only some of which are subalgebras of the Griess algebra.  We call an axial algebra with this Monster fusion law an axial algebra of \emph{Monster type}.
\begin{table}[!htb]
\setlength{\tabcolsep}{4pt}
\renewcommand{\arraystretch}{1.5}
\centering
\begin{tabular}{c||c|c|c|c}
 & $1$ & $0$ & $\frac{1}{4}$ & $\frac{1}{32}$ \\ \hline \hline
$1$ & $1$ &  & $\frac{1}{4}$ & $\frac{1}{32}$ \\ \hline
$0$ &  & $0$ &$\frac{1}{4}$ & $\frac{1}{32}$ \\ \hline
$\frac{1}{4}$ & $\frac{1}{4}$ & $\frac{1}{4}$ & $1, 0$ & $\frac{1}{32}$ \\ \hline
$\frac{1}{32}$ & $\frac{1}{32}$  & $\frac{1}{32}$ & $\frac{1}{32}$ & $1, 0, \frac{1}{4}$ 
\end{tabular}
\caption{Monster fusion law}\label{tab:monsterfusion}
\end{table}
Note that, since the Monster fusion law is $\mathbb{Z}_2$-graded, associated to each axis $a$, there is at most one non-trivial Miyamoto automorphism $\tau_a$ and it has order $2$.

Sakuma's theorem classifies the $2$-generated axial algebras of Monster type and shows that there are eight in total and they are labelled $2\A$, $2\B$, $3\A$, $3\C$, $4\A$, $4\B$, $5\A$ and $6\A$.  Here the number indicates the order of the product $\tau_a \tau_b$ in the Monster, where $a$ and $b$ generate the algebra.  In particular, the order of the product of any two Miyamoto involutions in an axial algebra of Monster type must have order at most $6$.  So, the Miyamoto group of an axial algebra of Monster type is a $6$-transposition group $(G, D)$, where each Miyamoto involution $\tau_a \in D$. 

In this note, we enumerate the axial algebras of Monster type whose Miyamoto group $G$ is minimal $3$-generated.  That is, where $G$ is $3$-generated, but any proper subgroup $H \lneqq G$ with $H = \la H \cap D \ra$ can be generated by two involutions in $D$.  Note that we must impose some additional restriction on a $3$-generated group $G$ to make the problem tractable.  Indeed, the Monster is $3$-generated and a $3$-generated $5$-transposition group $B(2,5):2$ contains an open case of the Burnside problem.  See \cite{3gen4trans} for an alternative restriction.

We note that this work was carried out independently of that of Mamontov, Staroletov and Whybrow in \cite{min3gen}.  However, both groups use the same list of minimal $3$-generated $6$-transposition groups which was originally given by Mamontov and Staroletov \cite[Theorem 2.2]{min3gen} and communicated to this author by Shpectorov at the Focused Workshop on Axial Algebras in May/June 2018.

We record our results here as they differ from theirs in the following ways:
\begin{itemize}
\item We do not assume that the algebra has a Frobenius form.  That is, a bilinear form which associates with the algebra product.  (Indeed, Mamontov, Staroletov and Whybrow assume a positive-definite Frobenius form.)
\item We do not assume primitivity i.e.\ that the $1$-eigenspace of every axis is $1$-dimensional.
\item We complete some cases which they could not (in addition to all those that they can).  Namely, we obtain two new non-trivial algebras and show a further seven do not lead to non-trivial algebras.  Furthermore, for one case we obtain a larger algebra which is a cover of the example in \cite{min3gen}.
\end{itemize}

For the construction of the algebras, we use the {\sc magma} implementation \cite{ParAxlAlg} of the algorithm in \cite{axialconstruction}, rather than the algorithm in \cite{Maddycode} used by Mamontov, Staroletov and Whybrow.

\begin{theorem}
If $A$ is a $3$-generated axial algebra of Monster type whose Miyamoto group is minimal $3$-generated, then it is a quotient of some algebra in Table $\ref{tab:results}$.
\end{theorem}

Note however that there are still twelve cases outstanding out of a total of 161.  Excluding these uncompleted cases and the two marked algebras computed over $\mathbb{Q}[t]$, we have the following:

\begin{corollary}\label{primFrob}
All the completed algebras in Table $\ref{tab:results}$ are primitive and admit a positive semi-definite Frobenius form.
\end{corollary}

This supports the conjecture in \cite{axialconstruction} that every axial algebra of Monster type admits a Frobenius form.

\section{Configurations of axes}

Given a $3$-generated $6$-transposition group $G = \langle x,y,z \rangle$, where $D = x^G \cup y^G \cup z^G$, we describe briefly how to find all possible actions on a putative set of axes $X$, so that $(A, X)$ could be an axial algebra of Monster type with Miyamoto group $G$.  For complete details, see Section 4 of \cite{3gen4trans}.

Using the action on $X$, we may assume that $a,b,c$ are axes such that $\tau_a = x$, $\tau_b = y$ and $\tau_c = z$ and $X = a^G \cup b^G \cup c^G$.  (Note that we do not assume that these orbits are disjoint.)  We begin by finding the possible stabilisers $G_d$ of an axis $d$, for $d = a,b,c$.  Note first that:
\[
\langle \tau_d \rangle \leq G_d \leq C_G(\tau_d)
\]
\begin{lemma}\textup{\cite[Lemma 4.3]{3gen4trans}}
$\tau_d = 1$ if and only if $G_d = G$.
\end{lemma}
We may now assume that $\tau_d$ is not the identity.  To help identify $G_d$, we make the following definitions.

\begin{definition}
\begin{enumerate}
\item An axis $d$ is \emph{unique} if there does not exist another axis $e \in X$ with $d \neq e$ and $\tau_d = \tau_e$.
\item An axis $d$ is \emph{strong} if there does not exist another axis $e \in d^G$ with $d \neq e$ and $\tau_d = \tau_e$.
\end{enumerate}
\end{definition}

Clearly a unique axis is strong.  Strong axes have axis stabilisers which are as large as possible:

\begin{lemma}\textup{\cite[Lemma 4.8]{3gen4trans}}
The following are equivalent
\begin{enumerate}
\item $d$ is strong
\item $G_d = C_G(\tau_d)$
\item There is a natural $G$-equivariant bijection between $d^G$ and $\tau_d^G$.
\end{enumerate}
\end{lemma}

With this in mind, we use the two following lemmas to gain a lower bound for $G_d$.

\begin{lemma}\textup{\cite[Lemma 4.5]{3gen4trans}}
Let $d \in X$.  If there exists $e \in X$ such that $\tau_d \tau_e$ has order $5$, then $d$ is unique.
\end{lemma}

\begin{lemma}\textup{\cite[Lemma 4.9]{3gen4trans}}
Let $d,e \in X$.
\begin{enumerate}
\item If the order of $\tau_d \tau_e$ is $2$ and $e$ is strong, then $\tau_e \in G_d$.
\item If the order of $\tau_d \tau_e$ is $4$, then $(\tau_d \tau_e)^2 \in G_d$.
\end{enumerate}
\end{lemma}

For $G = \langle x,y,z \rangle$, we enumerate all possible triples $(G_a, G_b, G_c)$ of axis stabilisers up to conjugation.  This naturally gives us a triple of orbits $(a^G, b^G, c^G)$, but we must still determine whether two orbits are disjoint or not.  If two of the generators $\tau_d$ and $\tau_e$ are conjugate and $G_d$ is conjugate to $G_e$, then there are two cases: either $d^G$ and $e^G$ are equal, or disjoint.  Otherwise, there is just one case and $d^G$ and $e^G$ are disjoint.  Taking into account all these different choices we may enumerate all possible actions of $G$ on the putative set of axes $X = a^G \cup b^G \cup c^G$.

As described in \cite{axialconstruction}, given a set of axes $X$, we may now construct the admissible $\tau$-maps.  That is, those maps $\tau \colon X \to G$ which could be a map giving the Miyamoto involutions $\tau_x$ for each axis $x \in X$.  In particular, $\Ima(\tau) = G$ and $\tau_{xg} = (\tau_x)^{g}$ for all $g \in G$.  Finally, for a given $X$ and $\tau$, we determine the set of possible shapes, where a shape is an assignment of $2$-generated subalgebra to each orbit of pairs of distinct axes.  This is the information needed for input into the algorithm.  Not that at each stage here in determining the action, $\tau$-map and shape, we do this up to appropriate automorphisms.

\section{Results}

We present a summary of our results in Table \ref{tab:summary}, giving the number of shapes for each group action and  indicate how many collapse, how many give a non-trivial completion and how many could not be completed.  In Table \ref{tab:results}, we give details of all the non-trivial algebras and all the cases which could not be completed.  Note that we do not include the Norton-Sakuma algebras here.  The columns in Table \ref{tab:results} are

\begin{itemize}
\item Miyamoto group.
\item Axes. We give the size as the sum of orbit lengths.
\item Shape.  We list one subalgebra for each connected component of the shape graph.
\item Whether the example is a minimal $3$-generated algebra in the sense that it does not contain any $3$-generated axial subalgebras.  For the cases we could not complete, we indicate whether it could be minimal or not.
\item Dimension of the algebra.  A question mark indicates that our algorithm did not complete and a $0$ indicates that the algebra collapses.
\item The minimal $m$ for which $A$ is $m$-closed.  Recall that an axial algebra is $m$-closed if it is spanned by products of length at most $m$ in the axes.
\item Whether the algebra has a Frobenius form that is non-zero on the set of axes $X$.  If it is additionally positive definite or positive semi-definite, we mark this with a pos, or semi, respectively.
\end{itemize}

\begin{table}[!htb]
\centering
\begin{tabular}{ccC{1.8 cm}C{1.8 cm}C{1.8 cm}C{1.9 cm}}
\hline
Group & Axes & Number of Shapes & Collapsing shapes & Shapes giving a non-trivial algebra & Incomplete shapes\\
\hline
$1$ & $1+1+1$ & 4 & 0 & 4 & \\
$2^2$ & $1+2+2$ & 6 & 0 & 6 & \\
$2^2$ & $2+2+2$ & 2 & 0 & 2 & \\
$2^2$ & $2+2+2$ & 6 & 2 & 3 & 1\\
$S_3$ & $1+3$ & 4 & 1 & 3 & \\
$S_3$ & $1+3+3$ & 3 & 0 & 3 & \\
$S_3$ & $3+3+3$ & 1 & 0 & 1 & \\
$2^3$ & $2+2+4$ & 4 & 4 & 0 & \\
$2^3$ & $2+2+4$ & 18 & 12 & 4 & 2\\
$2^3$ & $2+4+4$ & 12 & 11 & 1 & \\
$2^3$ & $4+4+4$ & 20 & 16 & 1 & 3\\
$D_{10}$ & $1+5$ & 2 & 1 & 1 & \\
$D_{12}$ & $2+6$ & 1 & 0 & 1 & \\
$D_{12}$ & $6+6$ & 2 & 2 & 0 & \\
$D_{12}$ & $2+6+6$ & 1 & 1 & 0 & \\
$D_{12}$ & $6+6+6$ & 2 & 2 & 0 & \\
$3^2:2$ & $9$ & 5 & 1 & 2 & 2 \\
$3^2:2$ & $9+9$ & 1 & 0 & 1 &  \\
$3^2:2$ & $9+9+9$ & 1 & 0 & 1 & \\
$S_4$ & $6$ & 4 & 0 & 4 & \\
$S_4$ & $3+6$ & 8 & 1 & 7 & \\
$S_4$ & $12$ & 1 & 0 & 1 & \\
$S_4$ & $6+6$ & 6 & 4 & 1 & 1 \\
$S_4$ & $3+6+6$ & 12 & 7 & 4 & 1 \\
$S_4$ & $6+6+6$ & 20 & 18 & 1 & 1 \\
$S_4$ & $12+12$ & 1 & 1 & 0 & \\
$S_4$ & $12+12+12$ & 1 & 1 & 0 & \\
$5^2:2$ & $25$ & 1 & 0 & 0 & 1 \\
$3^2:S_3$ & $27$ & 1 & 1 & 0 & \\
$3^2:S_3$ & $9+27$ & 1 & 1 & 0 & \\
$3^2:S_3$ & $9+9+27$ & 1 & 1 & 0 & \\
$3^2:S_3$ & $27+27$ & 1 & 1 & 0 & \\
$3^2:S_3$ & $9+27+27$ & 1 & 1 & 0 & \\
$3^2:S_3$ & $27+27+27$ & 1 & 1 & 0 & \\
$A_5$ & $15$ & 4 & 0 & 4 & \\
$5^2:S_3$ & $15$ & 2 & 0 & 1 & 1 \\
\hline
\end{tabular}
\caption{Summary}\label{tab:summary}
\end{table}

We now comment on our results.  Firstly, we complete all the cases completed in \cite{min3gen} and some additional cases.  In particular, we find two new non-trivial axial algebras, show an additional seven cases collapse and find a larger algebra for another case (this contains the algebra in \cite{min3gen} as a quotient).  (We also have a different result in one other case which we believe to be a typo in \cite{min3gen}.\footnote{We find that $S_3$ on $1+3+3$ axes with shape $6\A \, (2\A)^2$ has dimension $8$ and not $11$.  Indeed, it is isomorphic to the $6\A$ Norton-Sakuma algebra.\label{typo}})

The two new algebras are the $23$-dimensional and $36$-dimensional algebras for $S_4$ on $3+6+6$ axes with shapes $6\A \, 4\A \, (2\A)^3$ and $6\A \, 4\A \, 2\A \, 2\B \, 2\A$, respectively.  The cases which collapse are $2^3$ on $2+2+4$ axes with shape $(4\A)^2 \, (2\A)^2$; $2^3$ on $4+4+4$ axes with shapes $(4\A)^3 \, (2\A)^3$, $(4\A)^3 \, (2\A)^2 \, 2\B$, $(4\A)^3 \, 2\A \, (2\B)^2$, $(4\A)^2 \, 4\B \, (2\A)^3$ and $(4\A)^2 \, 4\B \, 2\A \, 2\B \, 2\A$ and $S_4$ on $6+6+6$ axes with shape $6\A \, (2\A)^6$.  Interestingly, one of the two new algebras is $3$-closed, but the other is $2$-closed and other $3$-closed algebras were constructed in \cite{min3gen}, so the difference between the two algorithms is not due to this.

In addition, we find a $16$-dimensional algebra $A$ for $2^3$ on $2+4+4$ axes with shape $4\B \, 4\A \, (2\A)^2$, whereas in \cite{min3gen} this is listed as $13$-dimensional.  However, our algebra $A$ has a positive semi-definite Frobenius form and when you factor out by the $3$-dimensional radical of this form, we recover a $13$-dimensional algebra with the same shape.  We believe that this accounts the difference here and confirms that Mamontov, Staroletov and Whybrow do indeed assume a positive definite Frobenius form.

For the constructed non-trivial algebras, it is easy to check that they are all primitive and admit a Frobenius form, giving us Corollary \ref{primFrob}.

\begin{longtable}{ccccccc}
\hline
$G$ & axes & shape & minimal & dim & $m$ & form \\
\hline
$1$ & 1+1+1 & $(2\A)^3$ & yes & 6,9 & 2,3 & pos \\
$1$ & 1+1+1 & $(2\A)^2 \, 2\B$ & yes & 6 & 3 & pos  \\
$1$ & 1+1+1 & $2\A \, (2\B)^2$ & yes & 4 & 2 & pos  \\
$1$ & 1+1+1 & $(2\B)^3$ & yes & 3 & 1 & pos \\
&&&&&\\

$2^2$ & 1+2+2 & $4\A \, (2\A)^2$ & no & 14 & 3 & semi \\
$2^2$ & 1+2+2 & $4\A \, 2\A \, 2\B$ & no & 10 & 3 & pos \\
$2^2$ & 1+2+2 & $4\A \, (2\B)^2$ & no & 6 & 2 & pos \\
$2^2$ & 1+2+2 & $4\B \, (2\A)^2$ & yes & 5 & 1 & pos\\
$2^2$ & 1+2+2 & $4\B \, 2\A \, 2\B$ & no & 8 & 2 & pos \\
$2^2$ & 1+2+2 & $4\B \, (2\B)^2$ & no & 6 & 2 & pos \\
&&&&&\\

$2^2$ & 2+2+2 & $4\A$  & yes & 12\,\footnote{Assuming primitivity, there is a 12 dim example over $\mathbb{Q}[t]$.  When evaluated for a particular value of $t$, the example always admits a Frobenius form, but it can be positive definite, positive semi-definite, or not either depending upon the value of $t$.  This was first constructed by Whybrow in \cite{maddyinfinite}.}
 & 3 & yes \\
$2^2$ & 2+2+2 & $4\B$ & yes & 7 & 2 & pos \\
&&&&&\\

$2^2$ & 2+2+2 & $4\A \, (2\A)^2$ & no & 0,10,?\,\footnote{Assuming primitivity, there are three possible cases when the algebra is evaluated over $\mathbb{Q}[t]$.  For all but two values of $t$, the algebra collapses.  When $t = \frac{1}{128}$, the algebra is 10-dimensional and it does not complete if $t=-\frac{1}{128}$.  For more details, see \cite{axialvarieties}} & 0,2,? & -,pos,? \\
$2^2$ & 2+2+2 & $4\A \, (2\B)^2$ & no & 9 & 3 & pos \\
$2^2$ & 2+2+2 & $4\B \, (2\A)^2$ & no & 11 & 2 & pos \\
$2^2$ & 2+2+2 & $4\B \, 2\A \, 2\B$ & no & 8 & 2 & pos \\
&&&&&\\

$S_3$ & 1+3 & $3\A \, 2\A$ & yes & 8 & 2 & pos \\
$S_3$ & 1+3 & $3\A \, 2\B$ & yes & 5 & 2 & pos \\
$S_3$ & 1+3 & $3\C \, 2\B$ & yes & 4 & 1 & pos \\
&&&&&\\

$S_3$ & 1+3+3 & $6\A \, (2\A)^2$ & no & 8\,\textsuperscript{\ref{typo}} & 2 & pos \\
$S_3$ & 1+3+3 & $6\A \, 2\A \, 2\B$ & no & 13 & 3 & pos \\
$S_3$ & 1+3+3 & $6\A \, (2\B)^2$ & no & 9 & 2 & pos \\
&&&&&\\

$S_3$ & 3+3+3 & $6\A$ & no & 13 & 2 & pos \\
&&&&&\\

$2^3$ & 2+2+4 & $(4\A)^2 \, 2\A\, (2\B)^2$ & no & ? & & \\
$2^3$ & 2+2+4 & $(4\A)^2 \, (2\B)^3$ & no & 13 & 3 & pos \\
$2^3$ & 2+2+4 & $4\A \, 4\B \, 2\A\, 2\B\, 2\A$ & no & 15 & 3 & pos\\
$2^3$ & 2+2+4 & $4\A \, 4\B \, (2\B)^2\, 2\A$ & no & 12 & 2 & pos \\
$2^3$ & 2+2+4 & $(4\B)^2 \, (2\A)^2\, 2\B$ & no & ? & & \\
$2^3$ & 2+2+4 & $(4\B)^2 \, (2\B)^3$ & no & 10 & 2 & pos \\
&&&&&\\

$2^3$ & 2+4+4 & $4\B \, 4\A \, (2\A)^2$ & no & 16 & 2 & semi \\
&&&&&\\

$2^3$ & 4+4+4 & $(4\A)^3 \, (2\B)^3$ & no & ? & & \\
$2^3$ & 4+4+4 & $(4\A)^2 \, 4\B \, (2\A)^2 \, 2\B$ & no & ? & & \\
$2^3$ & 4+4+4 & $4\A \, (4\B)^2 \, 2\B \, (2\A)^2$ & no & ? & & \\
$2^3$ & 4+4+4 & $(4\B)^3 \, (2\B)^3$ & no & 15 & 2 & pos\\
&&&&&\\

$D_{10}$ & 1+5 & $5\A \, 2\B$ & yes & 7 & 2 & pos \\
&&&&&\\

$D_{12}$ & 2+6 & $6\A$ & yes & 10 & 2 & pos \\
&&&&&\\

$3^2:2$ & 9 & $(3\A)^4$ & yes & ? & & \\
$3^2:2$ & 9 & $(3\A)^2 \, (3\C)^2$ & yes & ? & & \\
$3^2:2$ & 9 & $3\A \, (3\C)^3$ & yes & 12 & 2 & pos \\
$3^2:2$ & 9 & $(3\C)^4$ & yes & 9 & 1 & pos \\
&&&&&\\

$3^2:2$ & 9+9 & $6\A$ & no & 31 & 2 & pos \\
&&&&&\\

$3^2:2$ & 9+9+9 & $6\A$ & no & 42 & 2 & pos \\
&&&&&\\

$S_4$ & 6 & $3\A \, 2\A$ & yes & 13 & 2 & pos  \\
$S_4$ & 6 & $3\A \, 2\B$ & yes & 13 & 3 & pos \\
$S_4$ & 6 & $3\C \, 2\A$ & yes & 9 & 2 & pos \\
$S_4$ & 6 & $3\C \, 2\B$ & yes & 6 & 1 & pos \\
&&&&&\\

$S_4$ & 3+6 & $4\A \, 3\A \, 2\A$ & no & 23 & 3 & pos  \\
$S_4$ & 3+6 & $4\A \, 3\A \, 2\B$ & no & 25 & 3 & pos \\
$S_4$ & 3+6 & $4\A \, 3\C \, 2\B$ & no & 12 & 2 & pos  \\
$S_4$ & 3+6 & $4\B \, 3\A \, 2\A$ & yes & 13 & 2 & pos  \\
$S_4$ & 3+6 & $4\B \, 3\A \, 2\B$ & no & 16 & 2 & pos  \\
$S_4$ & 3+6 & $4\B \, 3\C \, 2\A$ & yes & 9 & 1 & pos \\
$S_4$ & 3+6 & $4\B \, 3\C \, 2\B$ & no & 12 & 2 & pos \\
&&&&&\\

$S_4$ & 12 & $6\A$ & yes & 17 & 2 & pos \\
&&&&&\\

$S_4$ & 6+6 & $6\A \, (2\A)^2 \, 2\B$ & no & 20 & 2 & pos \\
$S_4$ & 6+6 & $6\A \, (2\B)^2 \, 2\A$ & no & ? & & \\
&&&&&\\

$S_4$ & 3+6+6 & $6\A \, 4\A \, (2\A)^3$ & no & 23 & 2 & pos  \\
$S_4$ & 3+6+6 & $6\A \, 4\A \, 2\A \, 2\B \, 2\A$ & no & 36 & 3 & pos \\
$S_4$ & 3+6+6 & $6\A \, 4\A \, 2\B \, 2\A \, 2\B$ & no & ? & & \\
$S_4$ & 3+6+6 & $6\A \, 4\B \, (2\A)^2 \, 2\B$ & no & 20 & 2 & pos \\
$S_4$ & 3+6+6 & $6\A \, 4\B \, (2\B)^3$ & no & 23 & 2 & pos \\
&&&&&\\

$S_4$ & 6+6+6 & $6\A \, 2\A \, 2\B \, (2\A)^2 \, (2\B)^2$ & no & 28 & 2 & pos \\
$S_4$ & 6+6+6 & $6\A \, 2\B \, 2\A \, (2\B)^2 \, (2\A)^2$ & no & ? & & \\
&&&&&\\

$5^2:2$ & 25 & $(5\A)^6$ & yes & ? & & \\
&&&&&\\

$A_5$ & 15 & $5\A \, 3\A \, 2\A$ & yes & 26 & 2 & pos \\
$A_5$ & 15 & $5\A \, 3\A \, 2\B$ & no & 46 & 3 & pos \\
$A_5$ & 15 & $5\A \, 3\C \, 2\A$ & yes & 20 & 2 & pos \\
$A_5$ & 15 & $5\A \, 3\C \, 2\B$ & no & 21 & 2 & pos \\
&&&&&\\

$5^2:S_3$ & 15 & $5\A \, 3\A$ & yes & ? & & \\
$5^2:S_3$ & 15 & $5\A \, 3\C$ & yes & 18 & 2 & pos \\
\hline
\caption{$3$-generated axial algebras of Monster type with minimal $3$-generated Miyamoto groups} \label{tab:results}
\end{longtable}

\end{document}